\documentclass[12pt]{article} 
\usepackage{amsmath}
\usepackage[dvips]{graphicx}
\usepackage{amsfonts}
\usepackage{amsopn}
\usepackage{amsthm}
\def\Rn{{\mathbb R}^{n+1}}
\def\R3{{\mathbb R}^3}

\def\calO{{\cal O}}

\def\Sn{{\mathbb S}^n}
\def\R{{\mathcal R}}

\newcommand{\ep}{\epsilon}
 
\newcommand{\p}{\rho}

\newcommand{\n}{\nabla}

\newcommand{\ox}{\bar {x}}

\newcommand{\ow}{\bar {w}}

\renewcommand{\o}{\omega}

\renewcommand{\d}{\delta}
\newcommand{\D}{\Delta}

\newcommand{\ra}{\rangle}

\newcommand{\la}{\langle}

\newtheorem{theorem}{Theorem}

\newtheorem{lemma}[theorem]{Lemma}

\newtheorem{remark}[theorem]{Remark}
\newtheorem*{remark*}{Remark}
\title {A characterization of revolution quadrics by a system of partial differential equations}
\author{Vladimir I. Oliker\thanks{The research of the author
was partially supported by National Science Foundation grant DMS-04-05622.
}\\
\vspace{0.1in}\\
Department of Mathematics and Computer Science,\\
 Emory University, Atlanta, Georgia\\
oliker@mathcs.emory.edu}
\date{}
\begin{document}
\maketitle
\pagenumbering{arabic}
\setcounter{section}{0}
\setcounter{subsection}{0}
{\footnotesize
}
\begin{abstract}
It is shown  that  existence of a global solution to a particular nonlinear system of 
second order partial differential equations on a complete connected Riemannian manifold has topological and geometric implications and that in the domain of positivity of such solution its reciprocal  is the radial function of only one of the following rotationally symmetric hypersurfaces in $\Rn$:  paraboloid, ellipsoid, one sheet of a two-sheeted hyperboloid, and a hyperplane.
\end{abstract}
\maketitle
\pagenumbering{arabic}
\setcounter{section}{0}
\setcounter{subsection}{0}
{\footnotesize
}
\maketitle
\pagenumbering{arabic}
\setcounter{section}{0}
\setcounter{subsection}{0}
{\footnotesize
}

\section{Main result}\label{obata}
Let $M$ be a $C^{\infty}$ complete Riemannian manifold of dimension
$n \geq 2$ with metric $h$. Denote by $\n^2$ the Hessian matrix of the second covariant derivatives in the metric $h$. The well known 
theorem of Obata \cite{obata:62} states that
if the system
\begin{equation} \label{obata1}
 \n^2 w +  k^2wh = 0, ~k = const >0,
\end{equation}
admits a solution $w \not \equiv 0$
then $M$ is isometric to a sphere of radius $1/k$ 
in $\Rn$. 
There have been many important applications of Obata's theorem to problems of isometry and conformality  
of manifolds  with Euclidean spheres \cite{obata:71}, \cite{simon:77}, \cite{simon-lange:79}. On the other hand, some of the known characterizations of quadrics in Euclidean space use systems of third order partial differential equations (PDE's) defining second order spherical harmonics (see \cite{simon:07} and other references there), while
other characterizations are 
based on affine invariants  and affine isoperimetric inequalities \cite{li-simon:93}, \cite{Lut:86}, or on 
postulated relations between principal curvatures \cite{dillenetc:03}, \cite{kuhnel-book:06}, \cite{kuhnel-steller:05}, \cite{simon:07}. The proofs usually require substantial efforts.  By contrast, the characterization given in the following theorem uses a nonlinear second order system of PDE's, ultimately connected with first order spherical harmonics and it is obtained  by rather simple means 
akin to the original Obata's theorem.  In fact, the proof relies on Obata's theorem.

Recall, that a hypersurface $F$ in $\Rn$ is defined by its {\it radial}
function if $F$ is a graph of a positive function $\p$ over some domain $\o \subset \Sn$, where $\Sn$ is a unit sphere in $\Rn$ with the center at the origin of a Cartesian coordinate system. The position vector of $F$ is given by $\p(x)x,~x \in \o$. Here and below  $x$ is treated as a point in $\Sn$ and a unit vector in $\Rn$.
\begin{theorem}\label{thm}
Let $M$ be a $C^{\infty}$ complete connected Riemannian manifold of 
dimension $n \geq 2$ with metric $h$. 
 Suppose there exists a function $w \in C^2
(M), $ satisfying the system of PDE's
\begin{equation}\label{eq1}
2w\n^2w +(w^2-|\n w|^2) h = (c^2 -1)h ~\mbox{on}~ M,
\end{equation}
where $c = const$ and $w^2\not \equiv c^2-1$. 
Then $M$ is isometric to a unit sphere in $\Rn$.  Fix one such sphere $\Sn$ in $\Rn$  and chose a Cartesian coordinate system in $\Rn$ with the origin $\calO$ at the center of $\Sn$. Obviously, the function $-w$ is also a solution of (\ref{eq1}) and below we always deal with the solution $w$ which is positive somewhere on $\Sn$. The following statements hold:

(a) $~$ if $c^2 > 1$ then 
$w(x) =  \sqrt{C^2+c^2-1} +C 
\la x, \xi\ra,~x,~\xi \in \Sn,$ and $C=const\neq 0$. The function $w>0$ everywhere on $\Sn$.    The function $\p(x) = 1/w(x),~x \in \Sn,$ is the radial function of
an ellipsoid of revolution with axis $\xi$ and one of its foci at $\calO$. Depending on the choice of the center of $\Sn$ in $\Rn$ and $\xi$, such ellipsoid is defined up to a translation in $\Rn$ and rotation about the focus at $\calO$.
Furthermore, the eccentricity and location of the other focus are determined by the parameters $c$ and $C$. 

(b) $~$ if $c^2 =1$ then $w(x) = |C|+C\la x, \xi\ra, ~x,~\xi \in \Sn,$ and $C=const\neq 0$.  The function $w$ vanishes only at $\bar{x} =-sign(C)\xi$. The function $\p(x) = 1/w(x), ~x \in \Sn\setminus \{\bar{x}\},$ is the radial function of
a paraboloid of revolution with axis $\xi$ and focus at $\calO$. Again, depending on the choice of the center of $\Sn$ and $\xi$ this paraboloid is defined up to a translation in $\Rn$ and rotation about its focus. Replacing $C$ by $\lambda C$, where $\lambda > 0$, corresponds to a homothety of the paraboloid with respect to its focus.

(c) $~$ if $c^2 < 1$ then
$w^{\pm}(x) = \pm \sqrt{C^2+c^2-1} + C\la x,\xi\ra, ~x,~\xi \in \Sn,$ for $C^2=const \geq 1-c^2$. If $C^2> 1-c^2$ then $\p^{\pm}(x)=1/w^{\pm}(x),~x \in \o^{\pm}:=\{x\in \Sn~|~w^{\pm} > 0\},$ are the radial functions of one sheet of a two-sheeted hyperboloid of revolution with axis $\xi$ and one of the foci at $\calO$. Depending on the choice of the center of $\Sn$ in $\Rn$ and $\xi$, such hyperboloid is defined up to a translation in $\Rn$ and rotation about $\calO$. In addition, the eccentricity and location of the other focus are determined by the parameters $c$ and $C$. 
When $C^2= 1-c^2$ the corresponding function $\p(x)=1/C\la x,\xi\ra, ~x \in \o :=\{x \in \Sn~|~C\la x,\xi\ra > 0\}$, is the radial function of a hyperplane with normal $\xi$.
\end{theorem}

\begin{remark}Let $F$ be one of the following quadrics of revolution in $\Rn$: a paraboloid, an ellipsoid, and
one sheet of a two-sheeted hyperboloid of revolution. Assume that
the origin   of a Cartesian coordinate system in $\Rn$ is positioned at one of the foci of $F$. Then the radial function 
of $F$ is 
\begin{equation}\label{quad}
\p(x) = \frac{f}{1-\ep\la x,\xi\ra},~ x\in \o \subset \Sn,
\end{equation}
where
$\xi$ is a unit vector in the direction of the axis of $F$, $\ep$ is the eccentricity, and $f$ is the focal parameter. If $f > 0$ and $\ep \in [0,1)$ then $\o = \Sn$ and $F$ is an ellipsoid; when $f > 0,~\ep =1$ and $\o = \Sn\setminus \{\xi\}$ then $F$ is a paraboloid; if $f \neq 0$ and $|\ep| > 1$,
then (\ref{quad}) defines one of the sheets of a two-sheeted hyperboloid of revolution (over a sub-domain of $\Sn$ where $\p(x) > 0$). 

If $f \neq 0$ and $\p(x) = f/\la x,\xi\ra,~x \in \o :=\{x \in \Sn~|~f\la x,\xi \ra > 0\}$ then $F$ is a hyperplane. Note that in all cases the reciprocal $1/\p$ is defined on the entire $\Sn$ and satisfies (\ref{eq1}) with an appropriate choice of $c$.

It is not accidental that a one-sheeted hyperboloid of revolution is not among the hypersurfaces listed above. The reason is that it can not be represented
in the form (\ref{quad}) because such representation is possible only if the quadric of revolution has also a directrix hyperplane perpendicular to the axis of revolution. 
\end{remark}
\begin{remark} The system (\ref{eq1}) is motivated by the reflector problem and its generalizations \cite{OW}, \cite{refl_geom1}, \cite{oliker-refl:08}. Also, M. Gursky \cite{gursky:08} pointed out to the author that on a unit sphere $\Sn$ the left hand side of (\ref{eq1}) can be transformed into
an expression defining the Schouten tensor \cite{gursky-via:07}. More precisely, assume $w$ in (\ref{eq1}) is positive    and put $w=e^{u}$. Then $u$ satisfies the system
\begin{equation*}
e^{2u} \left [\n^2 u + \n u \otimes \n u + \frac{1-|\n u|^2}{2}h\right ] = \frac{c^2-1}{2}h.
\end{equation*}
The $(0,2)$-tensor in the square brackets is the Schouten tensor of the metric $e^{-2u}h$ on $\Sn$ with $h$ being the standard metric induced from $\Rn$.
\end{remark}
\begin{remark}\label{remark3}
In order to include characterizations of spheres of radius $\neq 1$ the statement in the first paragraph of the Theorem \ref{thm} should be modified as follows. Let $M$ be a complete connected Riemannian manifold of 
dimension $n \geq 2$ with metric $h$. Suppose there exists a function $w \in C^2
(M), $ satisfying the system
\begin{equation*}\label{eq-k}
2w\n^2w +(k^2w^2-|\n w|^2) h = (c^2 -1)h ~\mbox{on}~ M,
\end{equation*}
where $k=const > 0, ~c = const$ and $k^2w^2\not \equiv c^2-1$. 
Then $M$ is isometric to a sphere of radius $1/k$ in $\Rn$.
\end{remark}
\section{Proofs}\label{proof}
We will need the following 
\begin{lemma}\label{S0} Let $M$ be a complete connected Riemannian manifold of 
dimension $n \geq 2$ with metric $h$ and $A=const$. Suppose there exists a function $w \in C^2
(M)$, $w^2 \not \equiv A$, satisfying 
\begin{equation}\label{eq0}
2w\n^2w +(w^2-|\n w|^2 - A)h = 0 ~\mbox{on}~ M.
\end{equation}
Let $M_0=\{x \in M~|~w(x) =0\}$. (The case when 
$M_0=\emptyset$ is not excluded.) On the set $M\setminus M_0$
define a function $S$ by the equation
\begin{equation}\label{S}
w^2+|\n w|^2 + A = 2Sw.
\end{equation}
Then $S \equiv const$ on $M\setminus M_0$ and can be extended by the same constant to the entire $M$.
In addition, the function $w$ satisfies
the system
\begin{equation}\label{eq2}
\n^2w +(w-S)h = 0 ~\mbox{on}~ M.
\end{equation}
\end{lemma}
\begin{proof} Clearly, 
$S\in C^1(M\setminus M_0)$. Denote by $\n_i,~\n_{ij}$ the operators of the first and second covariant derivatives in some local coordinates.
Differentiating (\ref{S}), multiplying by $w$, and, taking into account (\ref{eq0}), we obtain after some manipulations
\[2w \n_i(Sw) = w^s[2w\n_{si} w + (w^2 - |\n w|^2 - A)h_{si}] + w^s[w^2 +|\n w|^2 + A]h_{si}= 2Sw\n_i w,\]
where $i=1,...,n, ~w^s = h^{sk}\n_k w, ~[h^{sk}]= [h_{sk}]^{-1}$ and the summation convention over repeated indices is in effect.
This implies $w^2\n_iS = 0, i=1,...,n,$ and we conclude that $S=const$ on $M\setminus M_0$.
Since $w\in C^2(M)$, we can extend $S$ to the entire $M$ by setting it equal to the same value
as on $M\setminus M_0$. Then the equation (\ref{S}) is satisfied everywhere on
$M$ with the same constant $S$.

To prove (\ref{eq2}) we calculate, using (\ref{S}),
\[
2w[\n^2w +(w-S)h] = 2w\n^2w +(w^2-|\n w|^2 - A)h =0. 
\]
This implies that 
\begin{equation}\label{eq3}
\n^2w +(w-S)h = 0 ~\mbox{on}~ M\setminus M_0.
\end{equation}
Let us show that the set $M_0$ has no interior points. It suffices to consider the case when $M_0\neq \emptyset$. Then the set $M\setminus M_0$ is open. Suppose, on the contrary, that
there exists an interior point $\ox \in M_0$. Take a geodesic $l(t)$ originating
at $\ox$ and parametrized by arc length. Then $w(l(t))=0$ for $|t| < \d$ for some $\d > 0$. This is true for any geodesic originating at $\ox$. In particular, if we take a geodesic joining $\ox$ with any point in $M\setminus \ox$, we obtain a contradiction as $w(l(t))$ will have to remain identically zero on such geodesic and this is impossible. Thus, $M_0$ can not have interior points and for each $\ox \in M_0$ there exists a sequence of points in $M\setminus M_0$ converging to $\ox$. Therefore, (\ref{eq3}) and $C^2$ continuity of $w$ imply 
(\ref{eq2}).
\end{proof}


\begin{proofof} {\bf Theorem \ref{thm}.}
Set $A=c^2-1$ in the Lemma \ref{S0}. Then the function $w$ satisfies (\ref{eq2}) with $S$ defined by (\ref{S}) everywhere on $M$. Consequently, the function $\ow := w-S$ satisfies 
the system
\[\n^2\bar{w} +  \bar{w}h = 0 ~\mbox{on}~ M.\]
If $w \equiv S$ then it follows from (\ref{S}) that $w^2\equiv c^2-1$. This contradicts one of the assumptions in the theorem. Therefore, $\ow \not \equiv 0$. 
By Obata's 
theorem $M$ is isometric to a sphere of unit radius in $\Rn$.
Taking the trace of this system we obtain
$$\D \bar{w} + n\bar{w} =0~\mbox{on}~ M.$$
It is well known that the only solutions of this equation on a unit sphere are  the first order spherical harmonics.
Fix one such sphere $\Sn$ in $\Rn$ and chose the Cartesian coordinate system in $\Rn$ with the origin at the center of $\Sn$.
Then $\bar{w}=C\la x,\xi\ra,$ where $C =const \neq 0, ~x\in \Sn,$ and $\xi$ is an arbitrary but fixed point of $\Sn$. Thus $w = S+C\la x,\xi\ra$. It is easy to check that $\bar{w}^2 + |\n \bar{w}|^2 = C^2$.
Then $(w-S)^2 + |\n w|^2 = C^2$ and using (\ref{S}) we obtain $S^2 = C^2+c^2-1$. 

Now we consider each of the cases in the statement of the theorem individually.

In case (a) $c^2-1 > 0$. By (\ref{S}) $Sw > 0$. Since it is assumed that $w > 0$, we also have $S > 0$ and then $S=\sqrt{C^2 + c^2 -1}$. Thus, any solution of (\ref{eq1}) is of the form $w = \sqrt{C^2 + c^2 -1} +C\la x, \xi\ra$. Put 
$$f= \frac{1}{\sqrt{C^2+c^2-1}},~~\ep= \frac{|C|}{\sqrt{C^2+c^2-1}}~~\mbox{and}~~\p(x) = \frac{f}{1\pm \ep\la x,\xi \ra}, ~x \in \Sn,$$
where $+$ is taken if $C> 0$ and $-$ if $C < 0$. Clearly, $\p$ is the radial function of an ellipsoid of revolution with axis $\xi$, focal parameter $f$
and eccentricity $\ep$. The remaining statements in (a)
follow from the construction. 

In case (b) $c^2-1=0$ and then, by one of the assumptions in the theorem $w \not \equiv 0$. Then $S \neq 0$  by (\ref{S}) and
$Sw \geq 0$. The choice of $w \geq 0$ on $\Sn$ implies $S >0$. Then $S=|C|$ and  any nonnegative solution of (\ref{obata1}) must be of the form $w(x) = |C| +C\la x, \xi\ra, ~x, ~\xi \in \Sn$. Clearly, for a fixed $\xi$ and $C\neq 0$ the function $w$ vanishes only at $x=-sign(C)\xi$.

Consider the
function 
$$\p(x) = \frac{1/|C|}{1+sign(C)\la x, \xi\ra},~~x \in \left\{\begin{array}{lr}
\Sn \setminus \{-\xi\}& \mbox{if} ~C > 0\\
\Sn \setminus \{\xi\} & \mbox{if}~C < 0. 
\end{array} \right.$$
In either case, this is the equation of radial function of a paraboloid  of revolution with axis $-\xi$ or $\xi$ and focal parameter $1/|C|$.
 Obviously, changing $\xi\in \Sn$ corresponds to rotation of such paraboloid
about its focus,
while multiplying $C$ by a positive constant corresponds to rescaling.

Finally, we consider the case (c). In this case $c^2-1 < 0$ and in order for
$S$ to be defined we require $C^2\geq 1-c^2$. Then $S=\pm \sqrt{C^2+c^2-1}$
and we have two representations $w^{\pm}(x) = \pm \sqrt{C^2+c^2-1} + C\la x,\xi\ra, ~x \in \Sn$.
Assuming first that $C^2+c^2-1 >0$ we put
\begin{eqnarray*}f= \pm\frac{1}{\sqrt{C^2+c^2-1}},~~\ep= \pm\frac{C}{\sqrt{C^2+c^2-1}}~~\mbox{and}
~~\p(x) = \frac{f}{1+\ep\la x,\xi \ra},\\
\mbox{where}~~~x \in\o:= \{x\in \Sn~|~w^{\pm} > 0\}.
\end{eqnarray*}
Note that $|\ep| > 1$ and $\p$ defines one sheet of a two-sheeted hyperboloid
of revolution.

It is easy to see that when $C^2+c^2-1 =0$ the hypersurfaces 
\[\p(x)=\frac{1}{C\la x,\xi\ra}, ~~x \in \o^{\pm}:=\left 
\{\begin{array}{lr}
\{x\in \Sn ~|~ \la x,\xi\ra > 0~\mbox{if}~C =\sqrt{1-c^2}\}\\
\{x\in \Sn ~|~ \la x,\xi\ra< 0~\mbox{if}~C =-\sqrt{1-c^2}\}
\end{array}
\right . 
\] 
are hyperplanes.
The remaining statements in case (c) are obvious. 
\end{proofof}

The {\bf proof of the claim in Remark \ref{remark3}} is based on a modified version of the Lemma  \ref{S0} and arguments similar to those in the first paragraph of the proof of  Theorem \ref{thm}. In the statement of the Lemma \ref{S0} the equations (\ref{eq0}), (\ref{S}) and (\ref{eq2}) should be replaced, respectively, by
 \begin{equation}\label{eq0prime}
2w\n^2w +(k^2w^2-|\n w|^2 - A)h = 0 ~\mbox{on}~ M,
\end{equation}
\begin{equation}\label{Sprime}
k^2w^2+|\n w|^2 + A = 2Skw
\end{equation}
and 
\begin{equation}\label{eq2prime}
\n^2w +k^2(w-\frac{S}{k})h = 0.
\end{equation}
The required changes in the proof of the so modified Lemma \ref{S0} are obvious.

\bibliographystyle{amsplain} 

\end{document}